\documentclass
{amsart}

\usepackage{amsmath,amsthm,amsfonts,amssymb,mathtools}
\usepackage{mathrsfs}
\makeatletter
\@namedef{subjclassname@2020}{\textup{2020} Mathematics Subject Classification}
\makeatother
\usepackage
{hyperref}
\usepackage{capt-of}
\usepackage{relsize}

\newcommand{\Be}{\begin{equation}}
\newcommand{\Ee}{\end{equation}}
\newcommand{\Bea}{\begin{eqnarray}}
\newcommand{\Eea}{\end{eqnarray}}
\newcommand{\Bel}{\begin{align}}
\newcommand{\Eel}{\end{align}}
\newcommand{\Beas}{\begin{eqnarray*}}
\newcommand{\Eeas}{\end{eqnarray*}}
\newcommand{\Benu}{\begin{enumerate}}
\newcommand{\Eenu}{\end{enumerate}}
\newcommand{\Bi}{\begin{itemize}}
\newcommand{\Ei}{\end{itemize}}
\newcommand\supp{\operatorname{supp}}

\def\R{{\mathbb R}}

\theoremstyle{plain}
\newtheorem{thm}{Theorem}[section]
\newtheorem{cor}[thm]{Corollary}
\newtheorem{lem}[thm]{Lemma}
\newtheorem{prop}[thm]{Proposition}

\theoremstyle{remark}
  
\theoremstyle{definition}

\numberwithin{equation}{section}

\newcommand{\RNum}[1]{\uppercase\expandafter{\romannumeral #1\relax}}

\usepackage{wrapfig}
\usepackage{tikz}

\newcommand{\ft}{{\mathfrak t}}

\begin{document}

\makeatletter
\@namedef{subjclassname@2020}{\textup{2020} Mathematics Subject Classification}
\makeatother
\subjclass[2020]{Primary 42B25,  Secondary 35S30}
\keywords{spherical maximal function, multi-parametric}

\author[Juyoung Lee]{Juyoung Lee}
\author[Sanghyuk Lee]{Sanghyuk Lee}
\author[Sewook Oh]{Sewook Oh}

\address{RIM, Seoul National University, Seoul 08826, Republic of Korea}
\email{ljy219@snu.ac.kr}

\address{Department of Mathematical Sciences and RIM, Seoul National University, Seoul 08826, Republic of Korea}
\email{shklee@snu.ac.kr}

\address{Korea Institute for Advanced Study, 85 Hoegiro Dongdaemun-gu, Seoul 02455, Republic of Korea}
\email{sewookoh@kias.re.kr}

\title[The strong spherical maximal function]
{$L^p$ bound on the strong  
spherical 
\\ maximal function}


\begin{abstract} In this note we show that  the strong spherical maximal function in $\mathbb R^d$ is bounded on $L^p$ if $p>2(d+1)/(d-1)$ for $d\ge 3$. 
\end{abstract}

\maketitle

\section{Introduction}
$L^p$ boundedness of maximal functions  of different forms, depending  on  specific purposes, has been studied in connection to  various related problems (see, for example, \cite{S} and \cite{CF, NSW}).
One of the most fundamental examples is  the Hardy--Littlewood maximal function, which is given by  maximal average over balls. 
For $\mathfrak{t}=(t_1,t_2,\cdots,t_d)\in\R^d$ and $y\in\R^d$, let 
\[ y_\mathfrak{t}=(t_1y_1, t_2y_2, \cdots, t_dy_d), \]
and set $m_{\mathfrak{t}}(f)=\int_{B(0,1)}f(x_{\mathfrak{t}})dx.$ Here $B(x,t)$ denotes the open ball centered at $x$ with radius $t>0$.
A natural multi-parametric generalization of the Hardy--Littlewood maximal function is  \emph{the strong maximal function} 
\[ 
\sup_{\mathfrak t\in \mathbb R^d}|f\ast m_{\mathfrak t}(x)|,\]
which   basically corresponds to the maximal function given by averages over the rectangles whose sides are parallel to the coordinate axes. As can be seen easily by $L^p$ boundedness of the one dimensional Hardy--Littlewood maximal function, 
it is well known that $f\mapsto  \sup_{\mathfrak t} |f\ast m_{\mathfrak t}| $ is bounded on $L^p$ for $p>1$. However, no $L^1$ theory is available for the strong maximal function \cite[p. 84]{S} (also see \cite{Ba,F,FGG}). 

Maximal averages over lower dimensional sets  also have attracted a lot of interest. Especially, the maximal functions given by averages over hypersurfaces have been extensively studied since Stein's seminal work \cite{S2}.  
For $d\ge 2$, let $\sigma$ denote the normalized surface measure on the unit sphere $\mathbb{S}^{d-1}$. By  $\sigma_t$ we  denote the measure $\sigma_t(f)=\int_{\mathbb{S}^{d-1}}f(ty)d\sigma(y)$. 
It was shown by Stein  \cite{S2} ($d\geq 3$) and Bourgain \cite{B} $(d=2)$  that  the spherical maximal operator $f\mapsto 
 \sup_{t>0}|f\ast d\sigma_t|$ 
  is bounded on $L^p$ if and only if $p>d/(d-1)$.

In this note we consider a multi-parameter version of the spherical maximal function. 
  For $\mathfrak{t}\in\R^{d}$,  let  $\sigma_{\mathfrak{t}}$ be the measure  given by 
  \[ \sigma_{\mathfrak{t}}(f)=\int_{\mathbb{S}^{d-1}}f(y_{\mathfrak{t}})d\sigma(y), \]
  which is a normalized measure on  the ellipsoid  $\mathbb S_{\mathfrak t} =\lbrace y_\mathfrak{t} : y\in\mathbb{S}^{d-1} \rbrace$. 
We consider a multi-parameter variant of the spherical maximal function
\[ \mathcal{M}f(x)=\sup_{\mathfrak{t}\in\R^d}|f\ast \sigma_{\mathfrak{t}}(x)|, \]
which we call the {\it strong spherical maximal function} in analogue of the strong maximal function. 
When $d=2$, it was shown in the authors' recent work \cite{LLO} that 
 $\mathcal{M}$  is bounded on $L^p$ for $p>4$. The range was later extended  to the optimal range $p>3$ by Chen, Guo, and Yang \cite{CGY}. 
 
 The objective of this note is to extend the 2-dimensional result to  higher dimensions, $d\geq 3$. The following is our main result. 
 
\begin{thm}\label{main thm}
   Let $d\geq 3$. For  $p>{2(d+1)}/(d-1)$,   there is a constant $C>0$ such that
   $ \Vert \mathcal{M}f\Vert_{L^p}\leq C\Vert f\Vert_{L^p} .$
\end{thm}

The range of $p$ is far from being optimal.  A Knapp type example shows that  $\mathcal M$ fails to be bounded on $L^p$ if $p<(d+1)/(d-1)$.  It seems possible  that the range can be further improved by elaborations on recent developments, especially, in decoupling inequalities, but 
it is unlikely that such attempts give the maximal bound on a sharp range of $p$. 
A straightforward consequence is a differentiation theorem which strengthens the earlier result \cite{S2}: 
\[ f\ast \sigma_{\mathfrak{t}}\to f  \quad{a.e.} \] 
as $\mathfrak t\to 0$ whenever $f\in L^p$ for $p>{2(d+1)}/(d-1)$.  Besides, let $E$ be a set of measure zero and $\tilde{\mathfrak t}:\R^d\to \R_+$ be a function. Then, it follows from 
the $L^p$ maximal bound that 
$|E\cap (\mathbb S_{\tilde{\mathfrak t}(x)}+x)|_{d-1}=0$ for a.e. $x$ (\cite{Marstrand}), where $|\cdot|_{d-1}$ denotes the $(d-1)$ dimensional Hausdorff measure. 

Our proof of Theorem \ref{main thm} relies on local smoothing estimates (see Proposition \ref{smoothing estimate} and \ref{two parameter weak smoothing} below) for multi-parameter wave operators
\[ \mathcal{U}_{\pm}f(x,\mathfrak{t})=a(x,\mathfrak{t})\int e^{i \Phi_{\pm}(x,\mathfrak{t},\xi)}\widehat{f}(\xi)d\xi, \]
where $a\in C_c^{\infty}(B(0,2)\times [2^{-1},2^2]^d)$ and $ \Phi_{\pm}(x,\mathfrak{t},\xi)=x\cdot\xi\pm |\xi_\mathfrak{t}|$.  The operators  $\mathcal{U}_{\pm}f$ are related to $f\mapsto f\ast\sigma_{\mathfrak t}$  via asymptotic expansion of the Fourier transform of the measure $\sigma_{\mathfrak t}$ (e.g., see  \eqref{ave} and \eqref{ft-sigma} below). 
The idea behind our approach is simple.  In perspective of one parameter sharp local smoothing estimates (e.g., see \eqref{1para})  for the wave operator  (combined with the Sobolev imbedding)   we only need  a small  extra gain of regularity  to obtain $L^p$ bound  on the $d$-parametric maximal function $\mathcal M$.  We show that this is possible by exploiting additional smoothing associated with  integration  in multi-parameters.  In $\R^2$, such smoothing estimates were shown to be true up to a sharp order in a certain range of $p$  \cite{LLO, CGY}.  Though it is an interesting problem  to obtain multi-parametric  local smoothing up to sharp order, we do not intend to pursue the issue here.  Instead, for our purpose, it is sufficient to obtain the smoothing estimates in Proposition \ref{smoothing estimate} and \ref{two parameter weak smoothing}. 
To prove those estimates we use decoupling inequalities (Proposition \ref{decoupling} and \ref{two para decoupling}), which are in turn to be proved  via variable coefficient generalizations of 
decoupling inequalities for conic surfaces  given by certain quadratic forms (Proposition \ref{key est} and \ref{key est two para}).  For this, we rely on recent developments concerning decoupling inequalities  due to Beltran--Hickman--Sogge \cite{BHS} and Guo--Oh--Zhang--Zorin-Kranich \cite{GOZZ}.

\subsection*{Notation}
In addition to the typical notation $\,\widehat{}\,\,$ we also use $\mathcal F$ to denote the Fourier transform. We denote $A\lesssim B$ if there exists a constant $C$ depending on $p,d$ such that $A\le CB$.

\section{Preliminaries and reduction}
We begin by  stating a lemma that plays the role of the Sobolev embedding (cf. Lemma 2.1 in \cite{LLO}).
\begin{lem}
\label{sobo3}
Let $1\le p\leq \infty$. For $1\leq i\leq d$, let $J_i$ be closed intervals of length $l\in [1/2, 2]$. 
Let $\mathfrak R=\prod_{i=1}^d J_i\subset \mathbb R^d$ and $G\in \mathrm C^1( \mathfrak R)$. Then, there is a constant $C>0$, independent of $\mathfrak R$ and $G$,  such that 
\begin{align*}
 & \sup_{v\in \mathfrak R}\vert G(v) \vert \le C (\prod_{i=1}^d\lambda_i)^{\frac{1}{p}}\sum_{\beta\in\{0,1\}^d}\tilde{\lambda}^{\beta}
        \Vert \partial_v^\beta G\Vert_{L^p(\mathfrak R)} 
 \end{align*}
for any  $\lambda_1, \dots, \lambda_d \ge 1$. Here  $\beta=(\beta_1, \beta_2, \cdots,\beta_d)\in \mathbb Z^d$ denotes a multi-index and $\tilde{\lambda}=(\lambda_1^{-1},\lambda_2^{-1},\cdots,\lambda_d^{-1})$. 
\end{lem}

 The lemma is easy to prove  by following the proof of Lemma 2.1 in \cite{LLO}, so we omit its proof.

We now consider  a local maximal operator
\[ \mathcal{M}_{loc}f(x)=\sup_{\mathfrak{t}\in[1,2]^{d}}|f\ast \sigma_{\mathfrak{t}}(x)|. \]
The estimates in the next proposition are the key estimates for our proof of  Theorem \ref{main thm}.

\begin{prop}\label{local maximal prop}  Let $d\ge 3$  and $j\ge 0$. Suppose $\supp\widehat{f}\subset \mathbb{A}_j:=\{ \xi : 2^{j-1}\leq |\xi|\leq 2^{j+1} \}$. Then, if $p> {2(d+1)}/({d-1})$, 
for some $\delta_0>0$  we have
   \Be \label{loc}
    \Vert \mathcal{M}_{loc}f\Vert_{L^p}\lesssim 2^{-\delta_0 j}\Vert f\Vert_{L^p}.
   \Ee
   \end{prop}

This proposition is to be proved by using Proposition \ref{smoothing estimate} and \ref{two parameter weak smoothing} below.

\subsection{Local smoothing estimates} 
By  Bourgain--Demeter's  decoupling inequality \cite{BD} and  scaling it is easy to see (e.g., see the proof of Lemma \ref{multipara lemma}) that 
 \Be
 \label{1para}\Vert \mathcal{U}_{\pm}f\Vert_{L^p_{x,\mathfrak{t}}}\lesssim 
        2^{(\frac{d-1}{2}-\frac dp+\epsilon)j}\Vert f\Vert_{L^p}, \quad \epsilon>0 
   \Ee
for $p\ge 2(d+1)/(d-1)$ whenever $\supp \widehat f\subset \mathbb A_j$. To get the maximal estimate \eqref{loc}, we need a  slightly stronger smoothing estimate than \eqref{1para}. 
To do this, due to some technical reasons (see Lemma  \ref{key-near} and \ref{key-away}),  we make use of two different types of local smoothing estimates depending 
whether $\widehat f$ is supported 
 near or away the coordinate axes.

 For a  small constant $c>0$, we define
\[ \mathbb{A}_j^{near}(c)= \bigcup_{k=1}^d \{ \xi\in \mathbb{A}_j:   ||\xi|^{-1} \xi -e_k|< c  \},\]
where $e_k$ denotes the $k$-th standard unit vector. We also set 
\[ \mathbb{A}_j^{away}(c)=\mathbb{A}_j\setminus  \mathbb{A}_j^{near}(c). \]

\begin{prop}\label{smoothing estimate}
    Let $d\ge 3$, $j\geq 0$ and $2\leq p\leq \infty$. Suppose that $\supp\widehat{f}\subset \mathbb{A}_j^{near}(c)$. Then, if $c$ is small enough,  for any $\epsilon>0$ we have
    \Be 
    \label{near}
    \Vert \mathcal{U}_{\pm}f\Vert_{L^p_{x,\mathfrak{t}}}\lesssim 
          2^{(\gamma(p)+\epsilon)j}\Vert f\Vert_{L^p},
      \Ee
   where
    \Be 
    \label{gamma}
    \gamma(p)=\begin{cases}
        \frac{d-1}{4}-\frac{d-1}{2p}, & 2\leq p\leq 6,
               \\[2pt]
        \frac{d-1}{2}-\frac{2d-2}{p}, & 6< p\leq \infty.\\
    \end{cases} 
    \Ee
\end{prop}

When $\supp \widehat f\subset \mathbb{A}_j^{away}(c)$, we have  estimates with  weaker smoothing but those are still better than the estimate   \eqref{1para}.

\begin{prop}\label{two parameter weak smoothing}
    Let $d\geq 3$, $j\geq 0$, and $2\leq p\leq \infty$. Suppose that $\supp\widehat{f}\subset \mathbb{A}_j^{away}(c)$ for a constant $c>0$. Then, for any $\epsilon>0$, we have
  \[
     \Vert \mathcal{U}_{\pm}f\Vert_{L^p_{x,\mathfrak{t}}}\lesssim 2^{(\alpha(p)+\epsilon)j}\Vert f\Vert_{L^p}, 
     \]
    where
    \Be \label{alpha-p} \alpha(p)=\begin{cases}
        \frac{d-1}{2}(\frac{1}{2}-\frac{1}{p}), & 2\leq p\leq \frac{2d}{d-2},
        \\[2pt]
        \frac{2d-3}{4}-\frac{2d-1}{2p}, & \frac{2d}{d-2}<p\leq 6,
        \\[2pt]
        \frac{d-1}{2}-\frac{d+1}{p}, & 6< p\leq \infty.\\
    \end{cases} \Ee
\end{prop}

The  estimates in Proposition \ref{two parameter weak smoothing} have smoothing of lower orders 
 than  those  in Proposition \ref{smoothing estimate}. In fact, to prove Proposition \ref{two parameter weak smoothing} we only 
 exploit local smoothing phenomena related to two parameters $t_j, t_k$, $1\le j<k\le d$ (see Section \ref{dddd}).

\subsection{Proof of Proposition \ref{local maximal prop}}
Once we have  Proposition \ref{smoothing estimate} and \ref{two parameter weak smoothing}, the proof of Proposition \ref{local maximal prop} is rather straightforward. 
Assuming  Proposition \ref{smoothing estimate} and \ref{two parameter weak smoothing} for the moment, we prove Proposition \ref{local maximal prop}. We begin by  recalling 
\Be
\label{ave}
\mathcal F (f\ast\sigma_{\mathfrak{t}})(\xi)=  
    \widehat{f}(\xi)\,    \widehat{\sigma}(\xi_{\mathfrak t}) 
\Ee
and 
the well known fact that
\Be \label{ft-sigma} \widehat{\sigma}(\xi)=e^{i|\xi|}a_+(\xi)+e^{-i|\xi|}a_-(\xi),  \quad |\xi|\ge 1 \Ee
and $a_\pm$ satisfy  $|\partial^l a_{\pm}(\xi)|\lesssim (1+|\xi|)^{-\frac{d-1}{2}-l}$ for any $l$.

    Let $p>{2(d+1)}/(d-1)$ and $\supp \widehat f\subset \mathbb A_j$. 
    Since $\ft\in [1, 2]^d$,   by \eqref{ave} and \eqref{ft-sigma} it is sufficient to consider  $2^{-(d-1)j/{2}} \mathcal{U}_{\pm}f$ instead of $f\ast\sigma_{\mathfrak{t}}$.  Note that
    \[ \partial_{t_i}\mathcal{U}_{\pm}f(x,\mathfrak{t})=\int e^{i(x\cdot\xi\pm|\xi_{\mathfrak{t}}|)} \Big(\partial_{t_i}a(x,\mathfrak{t})\pm\frac{it_i\xi_i^2}{|\xi_{\mathfrak{t}}|}a(x,\mathfrak{t})\Big)\widehat{f}(\xi)d\xi. \]
   By Lemma \ref{sobo3}, the desired estimate \eqref{loc} follows if we show that
   \[
        2^{(\frac{d}{p}-\frac{d-1}{2})j}\Vert \mathcal{U}_{\pm}f\Vert_{L^p_{x,\ft}}\lesssim 2^{-\epsilon_0j}\Vert f\Vert_{L^p}
\]
    for some $\epsilon_0>0$.   Combining Proposition \ref{smoothing estimate} and \ref{two parameter weak smoothing}, we have
    $\Vert \mathcal{U}_{\pm}f\Vert_{L^p_{x,\ft}}\lesssim 2^{(\max\{\alpha(p),\gamma(p)\} +\epsilon)j}\Vert f\Vert_{L^p}$ for  any $\epsilon>0$ if $2\le p\le \infty$. 
    Therefore, the estimate  follows since 
  \[
      \max \{\alpha(p),\gamma(p)\}  <\tfrac{d-1}{2}-\tfrac{d}{p}
\]
  for $p>{2(d+1)}/(d-1)$.  \qed

\subsection{Decoupling inequalities}  
We move on to the proofs of Proposition \ref{smoothing estimate} and Proposition \ref{two parameter weak smoothing}.
To prove the propositions we make use of  decoupling inequalities  for $ \mathcal{U}_{\pm}$. 

To state the inequalities, we need some notation. For $\kappa\in(0, 1)$, let $\{ \Theta^{\kappa}_m \}_{m=1}^N$ denote a collection of finitely overlapping open caps on $\mathbb{S}^{d-1}$ of diameter $L\in (2^{-1}\kappa, 2\kappa)$ such that 
$\bigcup_{m=1}^N  \Theta^{\kappa}_m=\mathbb S^{d-1}$.  Let $\{\zeta_m^\kappa\}_{m=1}^N$ be a partition of unity on $\mathbb{S}^{d-1}$
satisfying $\supp\zeta_m^{\kappa}\subset {\Theta}_{m}^{\kappa}$ for $1\le m\le N$ and $|\partial^l\zeta_m^{\kappa}|\lesssim \kappa^{-|l|}$ for any multi-index $l$. We denote 
\[ \mathfrak{S}(\kappa)=\{  \zeta_m^{\kappa} \}_{m=1}^N.\]
For each $\nu \in \mathfrak{S}(\kappa)$, set 
\[ \widehat{f_{\nu}}(\xi)=\widehat{f}(\xi)\nu (\xi/|\xi|).\]

Proposition \ref{smoothing estimate} is a consequence of the following decoupling inequality.

\begin{prop}\label{decoupling}
    Let $p\geq 6$ and $j\ge 0$. Suppose $\supp\widehat{f}\subset \mathbb{A}_j^{near}(c)$. Then, if $c>0$ is small enough, for any $\epsilon>0$ and $M>0$ we have
    \[ \Vert \mathcal{U}_{\pm}f\Vert_{L^p_{x,\mathfrak t}}\lesssim 2^{(\frac{d-1}{2}-\frac{2d-2}{p}+\epsilon)j}
    \Big(\sum_{\nu\in\mathfrak{S}(2^{-j/2})}\Vert \mathcal{U}_{\pm}f_{\nu}\Vert_{L^p_{x,\mathfrak{t}}}^p\Big)^{\frac{1}{p}}+2^{-Mj}\Vert f\Vert_{L^p}. \]
\end{prop}

Assuming this, we prove Proposition \ref{smoothing estimate}.

\begin{proof}[Proof of Proposition \ref{smoothing estimate}]
   The estimate \eqref{near} for $p=2$ is clear  by Plancherel's theorem. By interpolation, it suffices to show \eqref{near} for $p\geq 6$. 
 Let us begin by noting that  the inequality  
   \[ (\sum_{\nu\in\mathfrak{S}(2^{-j/2})}\Vert g_{\nu}\Vert_{L^p}^p)^{\frac{1}{p}}\lesssim \Vert g\Vert_{L^p}, \quad 2\le p\le \infty\]  
   holds when $\supp \widehat g\subset \mathbb A_j$. This can be shown by interpolation between the estimates for $p=2$ and $p=\infty$. Therefore, by Proposition \ref{decoupling} it is sufficient to show 
    \[ \Vert \mathcal{U}_{\pm}f_{\nu}\Vert_{L^p_{x,\mathfrak{t}}}\lesssim \Vert f_{\nu}\Vert_{L^p}, \quad 2\le p\le \infty \]
    for $\nu\in\mathfrak{S}(2^{-j/2})$.  Indeed, fix $\nu\in\mathfrak{S}(2^{-j/2})$ and $\mathfrak{t}\in [1,2]^d$. 
    Let  $\tilde{\mathfrak{t}}=(t_1^{-1},t_2^{-1},\cdots,t_d^{-1})$.  Changing  variables $\xi\to \xi_{\tilde{\mathfrak{t}}}$ and $x\to x_{\mathfrak{t}},$
    we see
    \[ \Vert \mathcal{U}_{\pm}f_{\nu}\Vert_{L^p_{x,\mathfrak{t}}(\R^d\times[1,2]^d)}\lesssim  \Big \Vert \int e^{i(x\cdot\xi\pm|\xi|)}\nu\big({\xi_{\tilde{\mathfrak{t}}}}/{|\xi_{\tilde{\mathfrak{t}}}|}\big)\mathcal Ff(\cdot_{\mathfrak{t}})(\xi)d\xi\Big\Vert_{L^p_{x,\mathfrak{t}}(\R^d\times[1,2]^d)} .\]
  The support of the function $\nu({\xi_{\tilde{\mathfrak{t}}}}/{|\xi_{\tilde{\mathfrak{t}}}|})\mathcal Ff(\cdot_{\mathfrak{t}})(\xi)$ is contained in a tube of width $\sim 2^{j/2}$ and length $\sim 2^j$. 
   By rotation, we may assume that $\nu({\xi_{\tilde{\mathfrak{t}}}}/{|\xi_{\tilde{\mathfrak{t}}}|})\mathcal Ff(\cdot_{\mathfrak{t}})(\xi)$ is supported in
    \[  R_c:=\{ \xi :  c^{-1}2^j\le  |\xi_1|\le c2^j, |\xi_i|\le c 2^{\frac{j}{2}}, 2\leq i\leq d \}\] 
    for some constant $c>0$. We choose  a smooth bump function $\tilde{\nu}$ adapted to  $R_c$. That is to say, $\tilde{\nu}=1$ on $R_c$, $\supp \tilde\nu\subset R_{2c}$, 
    and $\partial_1^\alpha(\partial_2\dots\partial_d)^\beta \tilde\nu\lesssim 2^{-\alpha j} 2^{-|\beta|j/2}$. 
   Using those properties of $\tilde \nu$, one can easily  see that 
   $ \Vert (e^{i(\pm|\xi|\mp \xi_1)}\tilde{\nu}(\xi))^{\vee}\Vert_{L^1}\lesssim 1$, equivalently, $ \Vert (e^{i(\pm|\xi|)}\tilde{\nu}(\xi))^{\vee}\Vert_{L^1}\lesssim 1. $ Thus,   by Young's convolution inequality and a simple change of variables, we obtain 
       \[
      \Vert \mathcal{U}_{\pm}f_{\nu}(\cdot, \mathfrak t)\Vert_{L^p_{x}} 
           \lesssim \Vert f_{\nu}\Vert_{L^p}.
    \]
    This gives the desired estimate by taking integration over $\ft \in [1,2]^d$. 
    \end{proof}

To prove Proposition \ref{two parameter weak smoothing}, we use the following decoupling inequality.
\begin{prop}\label{two para decoupling}
    Let $2\leq p\leq \infty$. Suppose $\supp\widehat{f}\subset \mathbb{A}_j^{away}(c)$ for a constant  $c>0$. 
     Then, for any $\epsilon>0$ and $M>0$ we have
    \[ \Vert \mathcal{U}_{\pm}f\Vert_{L^p_{x,\mathfrak{t}}}\lesssim 2^{(\alpha(p)+\epsilon)j}\Big(\sum_{\nu\in\mathfrak{S}(2^{-j/2})}\Vert \mathcal{U}_{\pm}f_{\nu}\Vert_{L^p_{x,\mathfrak{t}}}^p\Big)^{\frac{1}{p}}+2^{-Mj}\Vert f\Vert_{L^p}. \]
\end{prop}

One can prove  Proposition \ref{two parameter weak smoothing} using Proposition \ref{two para decoupling} in the same manner as Proposition \ref{smoothing estimate}, so we omit its proof.

\subsection{Decoupling inequalities for extension operators}
\label{dddd}
To prove Proposition  \ref{decoupling} and \ref{two para decoupling}, we need only to consider $\mathcal{U}_+$ since $\mathcal{U}_-$ can be handled in the same manner. By decomposition and  symmetry, we may assume that
\[ \supp\widehat{f}\subset \{ \xi \in \mathbb A_j: \xi_d>|\xi'|/\sqrt d\, \}, \]
where $\xi=(\xi',\xi_d)\in\R^{d-1}\times\R$. 
We write $\mathfrak t=(\mathfrak t', t_d)$ and fix $t_d\in [2^{-1},2]$.  For the proof of Proposition \ref{decoupling}  it is enough to show  the $\ell^p$-decoupling of the operator $\mathcal{U}_+ f(\cdot,  t_d)$  into  $\mathcal{U}_+ f_\nu(\cdot,  t_d)$  in $L^p_{x,\ft'}$ with  bound $C2^{(\frac{d-1}{2}-\frac{2d-2}{p}+\epsilon)j}$. 

To do  this, we make additional reduction. 
Changing variables $\xi\to 2^j\xi$ and $\xi\to (r\theta,r)$ $((\theta,r)\in \mathbb R^{d-1}\times \mathbb R)$, we write
\Be
\label{u+}
\mathcal{U}_+ f(x,\mathfrak t)=a(x,\mathfrak t)\iint e^{i2^jr(x'\cdot \theta+x_d+|(\theta,1)_{\mathfrak t}|)}\mathcal F[{f(2^{-j}\cdot)}](r\theta,r)r^{d-1}drd\theta. 
\Ee
 Performing harmless change of variables $(t^2_1, \dots, t_d^2)\to (\ft', s)$ and fixing $s  \in [2^{-2}, 2^4]$, we consider the operator
\[  \mathcal{E}^{s}_\lambda g(x,\mathfrak t')= \mathfrak a(x,\mathfrak t')\iint e^{i\lambda r(x'\cdot \theta+x_d+\Phi_s^{ne}(\ft', \theta))}g(\theta, r) drd\theta, \]
where $\mathfrak a\in C_c^{\infty}(B(0,2)\times [2^{-2},2^4]^{d-1})$  and
\Be 
\label{phi-s}
 \Phi_s^{ne}(\ft',\theta)=\textstyle \sqrt{s+\sum_{j=1}^{d-1}t_j\theta_j^2}. 
 \Ee
Here,  we use  the notation $ \Phi_s^{ne}$ to indicate that the phase function is associated with the part \emph{near} the coordinate axes. 

For $\kappa\in(0,1/(10d))$, let $\{ \mathfrak{q}_m^\kappa \}_{m=1}^N$ denote a collection of finitely overlapping open cubes of sidelength $L\in(2^{-1}\kappa,2\kappa)$ which cover $\{ \theta : |\theta|<\sqrt d\, \}$. 
Let $\{ \tilde{\zeta}_m^{\kappa}\}_{m=1}^N$ be a partition of unity subordinated to  $\{ \mathfrak{q}_m^\kappa \}_{m=1}^N$ such that  
$|\partial^l\tilde{\zeta}_m^\kappa|\lesssim \kappa^{-|l|}$ for multi-index $l$. We denote 
\[ \mathfrak{Q}(\kappa)=\{ \tilde{\zeta}_m^{\kappa} \}_{m=1}^N, \]
and for each $\mathfrak{q}\in\mathfrak{Q}(\kappa)$ we set
\[ g_\mathfrak{q} (\theta,r)=g(\theta,r)\mathfrak{q}(\theta). \]

Proposition \ref{decoupling} is a consequence of the next.

\begin{prop}\label{key est}
Let $p\geq 6$ and $s\in (2^{-2}, 2^4)$.  Suppose  $\supp_r g\subset [2^{-1}, 2]$ and $\supp_{\theta} g\subset B(0,\sqrt{c})$.
If $c$ is small enough,   for any $\epsilon>0$ and $M>0$ there is a constant $C$ such that 
    \[ \Vert \mathcal{E}_{\lambda}^s g\Vert_{L^p_{x,\mathfrak{t}'}}\le C \lambda^{\frac{d-1}{2}-\frac{2d-2}{p}+\epsilon}(\sum_{\mathfrak{q}\in\mathfrak{Q}(\lambda^{-1/2})}\Vert \mathcal{E}_{\lambda}^sg_\mathfrak{q} \Vert_{L^p_{x,\mathfrak{t}'}}^p)^{\frac{1}{p}}+\lambda^{-M}\Vert g\Vert_{L^p}. \]
\end{prop}

Set $\mathfrak t''=(t_1,\dots\,t_{d-2})$. To prove Proposition \ref{two para decoupling}, recall \eqref{u+}. 
 Changing variables $\mathfrak t\to ( t(\mathfrak t'', 1), ts)$,  
we may consider the operator $f\mapsto \mathcal U_+ f(\cdot, t(\mathfrak t'', 1), ts)$ instead of   $\mathcal U_+$. 
Thus,  fixing  $\mathfrak t''$, it is sufficient to obtain  $\ell^p$-decoupling of the consequent 2-parameter operator  
$f\mapsto \mathcal U_+ f(\cdot, t(\mathfrak t'', 1), ts)$ 
 in $L^p_{x,t,s}$ with bound $C2^{(\alpha(p)+\epsilon)j}$.    Recalling \eqref{u+} again, note that 
 $|(\theta,1)_{\mathfrak t}|$ is now replaced by $ t (|\theta_{(\mathfrak t'', 1)}|^2+ s^2)^{1/2}$. 
  By scaling in $\theta$  we may assume that $t_1=\cdots=t_{d-2}=1$. 
 Therefore, the matter is now reduced to showing an $\ell^p$-decoupling inequality for  the operator
\[  \tilde{ \mathcal{E}}_{\lambda}g(x,t,s)=\tilde{\mathfrak a}(x,t,s)\int e^{i\lambda r(x'\cdot\theta+x_d+\Phi^{aw}(t,s,\theta))}g(\theta,r)drd\theta. \]
where $\tilde {\mathfrak a}\in C_c^{\infty}(B(0,2)\times [2^{-2},2^4]^2)$ and 
\Be
\label{phi-aw} \Phi^{aw}(t,s,\theta):=t\sqrt{|\theta|^2+s^2}. \Ee
Similarly as  before,  the phase $ \Phi^{aw}$ is associated with the part \emph{away} from the axes. 
Unlike the operator $\mathcal{E}^{s}_\lambda$, the operator $\tilde{ \mathcal{E}}_{\lambda}$ enjoys  rotational symmetry in $\theta$. 
This is to be useful later for proving the desired decoupling inequality. 

 Proposition \ref{two para decoupling} follows from the next proposition.
\begin{prop}\label{key est two para}
    Let $2\leq p\leq \infty$ and $c>0$. Suppose  $\supp_r g\subset [2^{-1}, 2]$ and $\supp_{\theta} g\subset   B(0,\sqrt d) \setminus B(0,\sqrt{c})$.
Then, for any $\epsilon>0$ and $M>0$ there is a constant $C$ such that
    \[ \Vert\tilde{ \mathcal{E}}_{\lambda}g\Vert_{L^p_{x,t,s}}\leq C\lambda^{\alpha(p)+\epsilon}(\sum_{\mathfrak{q}\in \mathfrak{Q}(\lambda^{-1/2})}\Vert\tilde{ \mathcal{E}}_{\lambda}g_{\mathfrak{q}}\Vert_{L^p_{x,t,s}}^p)^{\frac{1}{p}}+\lambda^{-M}\Vert g\Vert_{L^p}. \]
\end{prop}

We prove Proposition \ref{key est} and \ref{key est two para} in Section \ref{decoupling proof section}.

\section{Proof of Theorem \ref{main thm}}
To prove the main theorem, by scaling it suffices to prove that the maximal function  
\[ \tilde{\mathcal{M}}f(x)=\sup_{\mathfrak{t}\in(0,1]^{d}}|f\ast \sigma_{\mathfrak{t}}(x)| \]
is bounded on $L^p$ for $p>2(d+1)/(d-1)$. 
Before beginning the proof, we consider a maximal operator with less parameters.
For $1\leq n\le d-1$, we set 
    \[ \mathfrak{M}_nf(x)=\sup_{\substack{t_1, \dots, t_n\in [1,2]}}|f\ast \sigma_{\tilde{\mathfrak{t}}}(x) | \]
    where $\tilde{\mathfrak{t}}=(t_1, t_2,\cdots,  t_{n-1},  t_n, t_n,\cdots, t_n)\in \R^d$. Making use of the sharp local smoothing estimate for the one parameter wave operator, i.e., \eqref{1para},  we obtain the following.
    \begin{lem}\label{multipara lemma}
     Let $d\ge 3$,   $1\leq n\le d-1$, and $2<p<\infty$. Then,  we have 
        \[ \Vert \mathfrak{M}_n f\Vert_{L^p}\lesssim 2^{-\delta j}\Vert f\Vert_{L^p} \]
        for some $\delta>0$ provided that $\supp \widehat f\subset \mathbb A_j$. 
       \end{lem}
    \begin{proof}
    Let $2<p<\infty$. As before, via the Littlewood-Paley decomposition we combine 
    \eqref{ave}, \eqref{ft-sigma} and   Lemma \ref{sobo3} 
          with  $G$ only depending on $\tilde{\mathfrak{t}}'\coloneqq (t_1,\dots, t_n)$. Thus, it suffices to show that
        \begin{equation}\label{routine est}
            \Vert \int e^{i(x\cdot\xi+ |\xi_{\tilde{\mathfrak{t}}}|)}\widehat{f}(\xi)d\xi\Vert_{L^p_{x,\tilde{\mathfrak{t}}'}(\R^d\times [1,2]^n)}\lesssim 2^{(\frac{d-1}{2}-\frac{n}{p}-\delta) j}\Vert f\Vert_{L^p}
        \end{equation}
        for  some $\delta>0$ provided that $\supp \widehat f\subset \mathbb A_j$.  We perform the change of variables  
        \[ (t_1,t_2,\cdots,t_n)\mapsto (t_1t_n,t_2t_n,\cdots, t_{n-1}t_n,t_n) .\]
        Note that $t_1, \dots, t_n\sim 1$. Fixing $t_1,\cdots,t_{n-1}$ and then changing variables in the frequency variables $ \xi_i\to  t_i^{-1} \xi_i$, $1\le i\le n-1$, we may apply the local smoothing estimate for the one parameter wave operator  $f\mapsto  e^{it\sqrt{-\Delta}} f$. 
           Consequently, we obtain  
      \[    \Vert \int e^{i(x\cdot\xi+ |\xi_{\tilde{\mathfrak{t}}}|)}\phi_j(|\xi|)\widehat{f}(\xi)d\xi\Vert_{L^p_{x,\tilde{\mathfrak{t}}'}(\R^d\times [1,2]^n)}\lesssim  2^{(\frac{d-1}{2}-\frac{d}{p}+\epsilon) j}\Vert f\Vert_{L^p}\]
       for any $\epsilon>0$ and $p\ge 2(d+1)/(d-1)$  whenever $\supp \widehat f\subset \mathbb A_j$. Since $n\le d-1$, this gives \eqref{routine est} for $p\ge 2(d+1)/(d-1)$.  When $2<p< 2(d+1)/(d-1)$, the estimate \eqref{routine est}  follows by 
       interpolation with the  estimate  $ \Vert \int e^{i(x\cdot\xi+ |\xi_{\tilde{\mathfrak{t}}}|)} \widehat{f}(\xi)d\xi\Vert_{L^2_{x,\tilde{\mathfrak{t}}'}(\R^d\times [1,2]^n)}\lesssim\Vert f\Vert_{L^2}$, which is 
       a consequence of  Plancherel's theorem. 
    \end{proof}

However, the $d$-parameter spherical maximal operator $\mathcal{M}$ can not be handled similarly  as above since the order of one parameter local smoothing is not enough to give the maximal estimate. 
To overcome this, we exploit  extra smoothing due to the multi-parameter averaging.

  Let  $\phi$ denote a smooth function such that $\supp \phi\subset (1-2^{-10},2+2^{-10})$ and $\sum_{j=-\infty}^{\infty}\phi(\tau/2^j)=1$ for $\tau>0$. 
For simplicity, we set  $\phi_j(\tau)=\phi(\tau/2^j)$,  
\[ \textstyle \phi_{<k}=\sum_{j<k}\phi_j, \quad \phi_{\geq k}= \phi- \phi_{<k}.\]
Let 
   \[Q_{\mathbf{k}}=\prod_{i=1}^d [2^{-k_i},2^{-k_i+1}].\]
    For $\omega=(\omega_1,\cdots,\omega_d)\in\{ 0,1\}^d$ and $\mathbf{k}=(k_1,k_2,\cdots,k_d)\in\mathbb N^d_0$,  define a projection operator
    \[ \mathcal F(\mathcal{P}_{\omega}^{\mathbf{k}}f)(\xi)= \prod_{i=1}^d\Big((1-\omega_i)\phi_{<k_i}(|\xi_i|)+\omega_i\phi_{\geq k_i}(|\xi_i|)\Big) \widehat{f}(\xi). \]
    Hence, $\sum_{\omega\in\{0,1\}^{d}} \mathcal{P}_\omega^{\mathbf{k}} f=f$ for each $\mathbf k\in \mathbb N^d_0$. 
    Since $ \tilde{\mathcal{M}}f(x)\le  
    \sup_{\mathbf{k}\in\mathbb N^d_0} \sup_{\mathfrak{t}\in Q_{\mathbf{k}}}|f\ast \sigma_{\mathfrak{t}}(x)|$, we have
    \begin{align*}
        \tilde{\mathcal{M}}f(x) \leq \sum_{\omega\in\{0,1\}^{d}} \sup_{\mathbf{k}\in\mathbb N^d_0}\sup_{\mathfrak{t}\in Q_{\mathbf{k}}}|\mathcal{P}_\omega^{\mathbf{k}}f\ast \sigma_{\mathfrak{t}}(x)| \eqqcolon \sum_{\omega\in\{0,1 \}^d }\mathcal{M}_\omega f(x).
    \end{align*}
Thus, the matter is reduced to showing that  
\Be \label{mw}
\|\mathcal{M}_\omega f\|_p\lesssim \|f\|_p
\Ee
 holds for every $ \omega \in\{0,1 \}^d$ if $p>2(d+1)/(d-1)$.

  For $\omega \in\{ 0,1\}^d$, let $|\omega |=\omega _1+\dots+ \omega _d$.  
We claim that \eqref{mw} holds for $p>2$ whenever $|\omega|\le d-1.$ 
 We prove this by induction on $|\omega |$.  Consider the case $|\omega |=0$ first. In such a case, it is easy to see that $\mathcal{P}_{\omega }^\mathbf{k}f=f\ast K$ with a kernel $K$ satisfying
    \[ |K(x)|\lesssim K^{\mathbf{k}}(x)\coloneqq \prod_{i=1}^d 2^{k_i}(1+2^{k_i}|x_i|)^{-N} \]
    for any $N>0$. Thus, we have 
    \[ |\mathcal{P}_{\omega }^{\mathbf{k}}f\ast \sigma_{\mathfrak{t}}(x)|\lesssim |f|\ast K^{\mathbf{k}}(x), \quad \ft\in Q_{\mathbf k}, \]
    which implies $\mathcal{M}_{\omega }f(x)\lesssim M_sf(x)$ where $M_s$ is the strong maximal operator in $\R^d$. 
    Therefore, we have   \eqref{mw}
    for $1<p\leq \infty$, in particular, $p>2$. 
    
    We now consider the case $1\leq |\omega| \leq d-1$. Let $1\le m\le d-1$ and assume that   \eqref{mw} holds for $p>2$ whenever  $|\omega |<m$. 
    We show that the same holds true for $|\omega|=m$.  
    Let $\omega \in\{0,1 \}^d$ such that $|\omega|=m$.     By symmetry, we may assume 
    \[ \omega=\omega _m:=(\,\overbrace{1,1,\cdots,1}^{m}\,,0,\cdots,0). \]
   Setting $\Omega_m= \{ \omega \in\{0,1\}^d : \omega _i=0, m+1\leq i\leq d \}\setminus\{ \omega _m\} $, we have
    \[ \mathcal F(\mathcal{P}_{\omega _m}^{\mathbf{k}}f)(\xi)=\Big(\prod_{i=m+1}^d\phi_{<k_i}(|\xi_i|)\Big) \widehat{f}(\xi)-\sum_{\omega '\in \Omega_m} \mathcal F(\mathcal{P}_{\omega '}^{\mathbf{k}}f)(\xi). \]
   
   For $z=(z_1, \dots, z_d)\in \mathbb R^d$, we write \[ z'=(z_1, \dots, z_m), \quad z''=(z_{m+1}, \dots, z_{d}).\]
   Let us set $K^{\mathbf k''}(x'')= \prod_{i=m+1}^d 2^{k_i}(1+2^{k_i}|x_i|)^{-N}$. Similarly with the previous case, we note that $|\mathcal F^{-1}  (\prod_{i=m+1}^d\phi_{<k_i}(|\cdot_i|) \widehat{f}\,)|\lesssim 
   \int |f(x', z'')|   K^{\mathbf k''}(x''-z'') dz''$.  Thus,  if  $\ft \subset Q^{\mathbf k} $, we have
   \[
        |\mathcal{P}_{\omega _m}^{\mathbf{k}}f\ast \sigma_{\mathfrak{t}}(x)|\lesssim \iint_{\mathbb{S}^{d-1}} |f(x'-y_{\mathfrak{t}}',z'')|d\sigma(y) K^{\mathbf k''}(x''-z'') dz''
        +\sum_{\omega\in \Omega_m}\mathcal{M}_{\omega}f(x). 
        \]
     Let us set 
         \[
         \tilde{\mathfrak{M}}_mf(x)=\sup_{\substack{t_1, \dots, t_m>0}}|\int_{\mathbb{S}^{d-1}}f(x'-y_{\mathfrak{t}}',x'')d\sigma(y)|. 
         \]
             Recalling the definition of $\mathcal{M}_{\omega}$, we obtain 
        \[ \mathcal M_{\omega_m}f(x)\lesssim M_s''(\tilde{\mathfrak{M}}_mf(x',\cdot))(x'') +\sum_{\omega\in \Omega_m}\mathcal{M}_{\omega}f(x),\]
        where $M_s''$ denotes the strong maximal function in $\mathbb R^{d-m}$. 
   Since $|\omega|\le m-1$ for $\omega\in \Omega_m$, by the triangle inequality and  the induction  assumption
  it follows that  $\|\sum_{\omega '\in \Omega_m}\mathcal{M}_{\omega '}f\|_p\lesssim \|f\|_p$  for $p>2$. 
  For the first term, we use the following lemma, which  is an easy consequence of Lemma \ref{multipara lemma}. 
     
      \begin{lem}\label{reduced multipara lemma}
       Let $d>2$ and $1\leq m<d$.   Then, $\Vert \tilde{\mathfrak M}_mf\Vert_{L^p}\lesssim \Vert f\Vert_{L^p}$
        holds for $p>2$.
    \end{lem}
    
    Therefore, Lemma \ref{reduced multipara lemma} and  $L^p$ boundedness of the strong maximal function  give 
    $ \Vert \mathcal{M}_{\omega _m}f\Vert_{L^p}\lesssim \Vert f\Vert_{L^p} $
   for $2<p\leq \infty$.  Consequently,  we obtain \eqref{mw} for $2<p\leq \infty$ provided that $|\omega|\le d-1$. 
   
    Finally, we consider  \eqref{mw} for the case $|\omega|=d$, i.e., $\omega=(1, \dots, 1)$.  We note that
    \[ \mathcal{M}_{\omega}f(x)\leq \sup_{\mathbf{k}\in\mathbb N^d_0}\sum_{\mathbf{s}\in\mathbb N^d_0}\sup_{\mathfrak{t}\in Q_{\mathbf{k}}}| (\mathcal{P}^{\mathbf{k}+\mathbf{s}}f)\ast \sigma_{\mathfrak{t}}(x)| \]
    where $\widehat{\mathcal{P}^{v}f}(\xi)=\widehat{f}(\xi)\prod_{i=0}^d\phi_{v_i}(|\xi_i|)$, $v\in \mathbb Z^d$. Minkowski's inequality gives
    \begin{equation}\label{w_d est}
        \Vert \mathcal{M}_{\omega}f\Vert_{L^p}\lesssim \sum_{\bf s\in\mathbb N^d_0}\big( \sum_{\mathbf{k}\in\mathbb N^d_0}\Vert \sup_{\mathfrak{t}\in Q_{\mathbf{k}}}| (\mathcal{P}^{\mathbf{k}+\mathbf{s}}f)\ast \sigma_{\mathfrak{t}}| \Vert_{L^p}^p \big)^{\frac{1}{p}}.
    \end{equation}
    By scaling it is clear that   the $L^p$ bound on  $f\mapsto  \sup_{\mathfrak{t}\in Q_{\mathbf{k}}}| (\mathcal{P}^{\mathbf{k}+\mathbf{s}}f)\ast \sigma_{\mathfrak{t}}|$
  equals that on $f\mapsto \sup_{\mathfrak{t}\in Q_{0}}| (\mathcal{P}^{\mathbf{s}}f)\ast \sigma_{\mathfrak{t}}|$. Thus, using  Proposition \ref{local maximal prop}, we have
    \[ \Vert \sup_{\mathfrak{t}\in Q_{\mathbf{k}}}| (\mathcal{P}^{\mathbf{k}+\mathbf{s}}f)\ast \sigma_{\mathfrak{t}}| \Vert_{L^p}\lesssim 2^{-\delta\max\{ s_i : 1\leq i\leq d \}}\Vert \mathcal{P}^{\mathbf{k}+\mathbf{s}}f\Vert_{L^p} \]
    for some $\delta>0$ when $ p>2(d+1)/(d-1)$.  Therefore, combining this and \eqref{w_d est} yields
    \[  \textstyle  \Vert \mathcal{M}_{\omega}f\Vert_{L^p}\lesssim \sup_{\mathbf s\in\mathbb N^d_0 } \big( \sum_{\mathbf{k}\in\mathbb N^d_0}\Vert \mathcal{P}^{\mathbf{k}+\mathbf{s}}f\Vert_{L^p}^p \big)^{\frac{1}{p}}.\]
  Using the inequality 
 $   \big( \sum_{\mathbf{k}\in\mathbb N^d_0}\Vert \mathcal{P}^{\mathbf{k}+\mathbf{s}}f\Vert_{L^p}^p \big)^{\frac{1}{p}}\lesssim \Vert f\Vert_{L^p}$,  which holds 
    for $2\leq p\leq \infty$ (for example, see \cite[Lemma 6.1]{TVV}),   we obtain
    \eqref{mw} 
     for $p> {2(d+1)}/(d-1)$. This completes the proof.

\begin{proof}[Proof of Lemma \ref{reduced multipara lemma}]
Lemma \ref{reduced multipara lemma} can be shown  in a similar way by using the inductive argument (in the proof of Theorem \ref{main thm}) and Lemma \ref{multipara lemma} which replaces Proposition \ref{local maximal prop}. 
So, we omit the detail.
\end{proof}

\section{Decoupling inequalities}
\label{decoupling proof section}
In this section, we prove Proposition \ref{key est} and \ref{key est two para}. 
The proofs of the propositions are rather standard (\cite{BHS}). So, we shall be brief in what follows. 
We first deduce  decoupling inequalities for the extension operators given by the associated conic surfaces 
from those inequalities for quadratic surfaces. Afterward, we extend the resulting decoupling  inequalities to 
variable coefficient forms to get the desired decoupling inequalities.  

\subsection{Decoupling inequality for quadratic surfaces}
We begin by recalling a theorem concerning decoupling inequalities for quadratic surfaces. 

For an $n$-tuple of real quadratic forms $Q(\theta)=(Q_1(\theta),\cdots,Q_n(\theta))$, define an extension operator
\[
  E_\lambda^{Q}h(x',z)=\int e^{i\lambda (x'\cdot \theta+z \cdot Q(\theta))} h(\theta)d\theta, \quad (x',z)\in \mathbb R^{d-1}\times \mathbb R^n.
 \]

 Let  $
 \mathcal N(Q)=|\{1\le d'\le d-1:\partial_{\theta_{d'}} Q_{n'}\not\equiv 0 \ for \ some \ 1\le n'\le n\}|,$
 which is  the number of variables that $Q$ depends on.   For $0\le d'\le d-1$ and $0\le n'\le n$, set
 \[
\mathfrak d_{d',n'}(Q):=\inf_{\substack{M\in \R^{(d-1)\times (d-1)}\\
rank(M)=d'}}\inf_{\substack{M'\in \R^{n\times n}\\
rank(M')=n'}}\mathcal N(M'\cdot(Q\circ M)).
\]

\begin{thm}
[{\cite[Theorem 1.1]{GOZZ}}] \label{GOZZ}
Let $2\le p\le \infty.$ 
Let $Q$ be an $n$-tuple of real quadratic forms and
 \[
\alpha(p,Q):=\max_{0\le d'\le d-1}\max_{0\le n'\le n}\big(\big(d'-\tfrac{\mathfrak d_{d',n'}(Q)}2\big)\big(\tfrac 12 -\tfrac 1p\big)-\tfrac{n-n'}{p}\big).
\]
 Suppose $\supp h\subset B(0,1)$. Then, we have 
   \[
  \|E_\lambda^{Q} h\|_{L^p_{x',z} (w)}\lesssim \lambda^{\alpha(p,Q)+\epsilon}
  (\sum_{\mathfrak{q}\in\mathfrak{Q}(\lambda^{-1/2})}\Vert  E^Q_{\lambda} ( h\mathfrak{q}) \Vert_{L^p_{x',z}(w)}^p)^{\frac{1}{p}}.\]  
 \end{thm}

Here $w$ is a rapidly decaying weight function, which is concentrated in a unit ball, i.e.,  $w(x',z)=(1+|(x',z)|)^{-N}$ for a large $N$. 
Computing $\alpha(p,Q)$ for each $Q$ does not seem to be a simple matter. 
For our purpose,  we combine Theorem \ref{GOZZ} and  stability of the decoupling exponent that was shown  in \cite{GOZZ}. 
\begin{thm}[{\cite[Theorem 5.2]{GOZZ}}]\label{stability}
Let $2\le p\le \infty.$ 
Let $Q, \tilde Q$ be $n$-tuples of real quadratic forms and $\alpha(p, Q)$ be given as in Theorem \ref{GOZZ}.  Then, there exists $\sigma=\sigma(Q)>0$ such that 
\[
  \| E_\lambda^{\tilde Q} h\|_{L^p_{x',z}(w)}\lesssim \lambda^{\alpha(p,Q)+\epsilon}
  (\sum_{\mathfrak{q}\in\mathfrak{Q}(\lambda^{-1/2})}\Vert  E^{\tilde Q}_{\lambda} (h \mathfrak{q}) \Vert_{L^p_{x',z}(w)}^p)^{\frac{1}{p}}\]
 holds for  $\supp h\subset B(0,1)$ whenever  $\sup_{\theta\in B(0,1)}|\tilde Q(\theta)-Q(\theta)|<\sigma(Q)$.
\end{thm}

We now  consider a class of $n$-tuples of real quadratic forms associated to a smooth function $\Phi$ defined on a subset of $\R^{n}\times \R^{d-1}$. 
By  $Q(\Phi, z_0, \theta_0)$
we denote the $n$-tuple of quadratic forms $Q=(Q_1, \dots, Q_n)$ such that 
\[
\partial_{\theta_{k}}\partial_{\theta_{l}}Q_j=\partial_{\theta_{k}}\partial_{\theta_{l}}\partial_{z_j}\Phi(z_\circ,\theta_\circ), \quad j=1, \dots, n
\]
for each $1\leq k,l\leq d-1$. Combing Theorem \ref{GOZZ} and \ref{stability}, we have the following two lemmas.

\begin{lem}\label{key-near} Let $2\le p\le \infty.$ 
Let $Q=Q(\Phi^{ne}_s, \theta_\circ, z_\circ)$ for  $(\theta_\circ, z_\circ) \in  B(0,1)\times [2^{-2},2^4]^{d-1} $ and 
$\gamma(p)$ be given by 
\eqref{gamma}.
Suppose $|\theta_0|<c$ for some  $c>0$ small enough,  then for any $\epsilon>0$  we have
    \[ \Vert {E}^Q_{\lambda}  h\Vert_{L^p_{x',z}(w)}\lesssim \lambda^{\gamma(p)+\epsilon}(\sum_{\mathfrak{q}\in\mathfrak{Q}(\lambda^{-1/2})}\Vert {E}^Q_{\lambda}(h \mathfrak{q}) \Vert_{L^p_{x',z}(w)}^p)^{\frac{1}{p}}
       \]
 with the implicit constant independent  of $Q$  when $\supp h\subset B(0,1)$.
\end{lem}

\begin{proof}
Recalling  \eqref{phi-s}, it is easy to see that the $j$-th component $Q_j$ of 
$Q=Q(\Phi^{ne}_s, \theta_\circ, z_\circ)$
is given by 
\[  Q_j(\theta)=(2\sqrt s)^{-1} \theta_j^2 +O( c|\theta|^2) \]
for $|\theta_\circ|<c$. Thus, the quadratic forms $Q(\Phi^{ne}_s, \theta_\circ, z_\circ)$ are contained in $O(c)$-neighborhood of 
the $(d-1)$-tuple quadratic form 
\[  Q_0(\theta)= (2\sqrt s)^{-1} (\theta_1^2, \dots,\theta_{d-1}^2). \]   
One can easily check $\alpha(p, Q_0)=\gamma(p)$ (for example, see \cite[Corollary 3.3]{GOZZ}). Combing this and Theorem \ref{stability}, we obtain the desired inequality. 
\end{proof}

\begin{lem}\label{key-away} 
  Let $2\leq p\leq \infty$ and $c>0$. 
Let $Q=Q(\Phi^{aw}, \theta_\circ, z_\circ)$ for  $\theta_0\in B(0,2\sqrt d) \setminus B(0,\sqrt{c})$, $z_\circ\in [2^{-2},2^4]^{2} $, and 
$\alpha(p)$ be given by \eqref{alpha-p}. Then, for any $\epsilon>0$  we have
    \[ \Vert {E}^Q_{\lambda}  h\Vert_{L^p_{x',z}(w)}\lesssim \lambda^{\alpha(p)+\epsilon}(\sum_{\mathfrak{q}\in\mathfrak{Q}(\lambda^{-1/2})}\Vert {E}^Q_{\lambda}(h \mathfrak{q}) \Vert_{L^p_{x',z}(w)}^p)^{\frac{1}{p}}\]
 with the implicit constant independent  of $Q$  when $\supp h \subset B(0,1)$.
\end{lem}

\begin{proof}
Recall \eqref{phi-aw}. 
By rotation symmetry of $\Phi^{aw}$ in $\theta$ we may assume 
that $\theta_\circ=a e_1$ for some $a\ge c$.  Observe that  $\partial_{\theta_n}\partial_{\theta_m}\partial_{t}\Phi^{aw}(t,s,ae_1),  \partial_{\theta_n}\partial_{\theta_m}\partial_{s}\Phi(t,s,ae_1)=0$ whenever $n\neq m$. Here $z_\circ=(t,s)$. Furthermore,  a computation gives
\[
Q_\circ(\theta):=Q(\Phi^{aw}, \theta_\circ, z_\circ)(\theta) 
=\Big(Q^{t,s}(\theta),\frac{ts}{a^2+s^2}\Big(\frac{2a^2\theta_1^2}{(a^2+s^2)^{3/2}}-Q^{t,s}(\theta)\Big)\Big)
\]
where 
\[\textstyle Q^{t,s}(\theta)={(a^2+s^2)^{-\frac 32}}s^2\theta_1^2+(a^2+s^2)^{-\frac 12}\sum_{j=2}^{d-1}\theta_j^2.\]
Since $s,t\in [2^{-2},2^4]$, by harmless linear change of variables in $z$ and scaling in $\theta$, we may, in fact,  assume that  
$Q_\circ(\theta)=( \theta_2^2+\cdot+  \theta_{d-1}^2  ,  a^2\theta_1^2)$.  Consequently,  since $a\neq0$, it  is not difficult to see that  $\mathfrak d_{d',0}(Q_\circ)=0$, $\mathfrak d_{d',1}(Q_\circ)=\delta_{d',d-1},$  
and $\mathfrak d_{d',2}(Q_\circ)=d'$.  Thus, we have 
\[
\alpha(p,Q_\circ)=\max\big\{\tfrac{d-1}2 \big(1-\tfrac 2p\big)-\tfrac 2p, \big(d-1 -\tfrac12\big)\big(\tfrac12-\tfrac 1p\big)-\tfrac 1{p},\tfrac{d-1}2(\tfrac 12-\tfrac 1p)\big\}.
\]
Note that the right hand side is equal to   $\alpha(p)$ given by \eqref{alpha-p}. Combining this and Theorem \ref{stability}, we obtain the desired inequality. 
\end{proof}

\subsection{Conical extensions}
For an $n$-tuple $Q$ of quadratic forms, we define an extension operator
\[ 
\mathcal C ^Q_{\lambda}g(x,z)=\iint e^{i\lambda r(x'\cdot \theta+x_d+z\cdot Q(\theta))}g(\theta,r)drd\theta. 
\]
Combining Lemma \ref{key-near} and \ref{key-away} and the standard argument (\cite{BD}) which allows one to deduce the decoupling inequality for conic surfaces, 
we have the following  Corollary \ref{key-nearc} and \ref{key-awayc}. The corollaries  can be shown by routine adaptation of the projection argument in \cite{BD} that 
makes it possible to obtain the decoupling inequality for $\mathcal C_\lambda^Q$ from that for $E_\lambda^Q$. 
The argument is even clearer here since we are dealing with quadratic forms.  
We state the next two corollaries without  proofs.\footnote{It is also possible to give an alternative proof making use of Fourier series expansion by which one can disregard unwanted nonlinear part of the phase function.}

\begin{cor}\label{key-nearc} Let $2\le p\le \infty.$ 
Let $Q=Q(\Phi^{ne}_s, \theta_\circ, z_\circ)$ for  $(\theta_\circ, z_\circ) \in  B(0,1)\times [2^{-2},2^4]^{d-1} $ and 
$\gamma(p)$ be given by \eqref{gamma}.
Suppose $|\theta_0|<c$ for some  $c>0$ small enough,  then for any $\epsilon>0$  we have
    \[ \Vert \mathcal C^Q_{\lambda}g\Vert_{L^p_{x,z}(w)}\lesssim \lambda^{\gamma(p)+\epsilon}(\sum_{\mathfrak{q}\in\mathfrak{Q}(\lambda^{-1/2})}\Vert {\mathcal C}^Q_{\lambda}g_\mathfrak{q} \Vert_{L^p_{x,z}(w)}^p)^{\frac{1}{p}} 
    \]
 with the implicit constant independent  of $Q$  when $\supp_\theta g\subset B(0,1)$.
\end{cor}

\begin{cor}\label{key-awayc} 
  Let $2\leq p\leq \infty$ and $c>0$. 
Let $Q=Q(\Phi^{aw}, \theta_\circ, z_\circ)$ for  $\theta_0\in B(0,2\sqrt d) \setminus B(0,\sqrt{c})$, $z_\circ\in [2^{-2},2^4]^{2} $, and 
$\alpha(p)$ be given by \eqref{alpha-p}. Then, for any $\epsilon>0$  we have
    \[ \Vert \mathcal C^Q_{\lambda} g\Vert_{L^p_{x,z}(w)}\lesssim \lambda^{\alpha(p)+\epsilon}(\sum_{\mathfrak{q}\in\mathfrak{Q}(\lambda^{-1/2})}\Vert \mathcal C^Q_{\lambda} g_\mathfrak{q} \Vert_{L^p_{x,z}(w)}^p)^{\frac{1}{p}}\]
 with the implicit constant independent  of $Q$  when $\supp_\theta g\subset B(0,1)$.
\end{cor}

 Once we have Corollary \ref{key-nearc} and \ref{key-awayc}, which give uniform decoupling inequalities for the conic surfaces associated with the quadratic forms,  one can prove Proposition \ref{key est} and \ref{key est two para} following the argument in 
\cite{BHS}.  In fact, the phase functions of the  operators $\mathcal{E}^{s}_\lambda$ and $\tilde{\mathcal E}_\lambda$ are close to linear phases. 
Consequently, the proofs of Proposition \ref{key est} and \ref{key est two para} become even simpler.  

We briefly mention how Corollary \ref{key-nearc} and \ref{key-awayc} can be used for the proof of Propositions.  The proofs of Proposition \ref{key est} and \ref{key est two para} are essentially identical.  We only consider the former. 
Let us set $\Phi(z,\theta)=   \Phi_s^{ne}(z,\theta)$. Fixing $z_\circ\in[2^{-1},2]^{n}$ and $\theta_\circ \in B(0,1)$,  set $Q=Q(\Phi, \theta_\circ, z_\circ)$. 
   Note that
\begin{align*}
   &\Phi(z,\theta)-\Phi(z_\circ,\theta)
   -\Phi_z'(z_\circ,\theta_\circ)(z-z_\circ)
   -(\theta-\theta_\circ)^\intercal\Phi_{\theta, z}''(z_\circ,\theta_\circ)(z-z_\circ) \\ 
    &=(z-z_\circ)\cdot Q(\theta-\theta_\circ)+O(|z-z_\circ||\theta-\theta_\circ|^3)+O(|z-z_\circ|^2).
\end{align*}
Disregarding some harmless terms, after 
appropriate localization in $z,\theta$ and scaling, we may think of $\Phi(z,\theta)$ as a perturbation of 
$(z-z_\circ)\cdot Q(\theta-\theta_\circ)$ with some negligible error. Thus, using the uniform decoupling inequality in Corollary \ref{key-nearc} and following 
the argument in \cite{BHS} (also see \cite{LLO}), 
one can prove Proposition \ref{key est}.  We omit the details.

\subsection*{Acknowledgement} 
This work was supported by the National Research Foundation (Republic of
Korea) grant no. 2022R1A4A1018904 (J. Lee and S. Lee) and KIAS individual grant SP089101 (S. Oh).

\end{document}